\documentclass[12pt]{article}
\usepackage{hyperref}
\usepackage{cmap}
\usepackage[cp1251]{inputenc}
\usepackage[english]{babel}
\usepackage[left=1in,right=1in,top=1in,bottom=1.5in]{geometry}
\usepackage{amssymb,amsmath, amsthm, amscd,ifthen}
\usepackage{graphicx}

\newcommand{\R}{{\mathbb R}}
\newcommand{\N}{{\mathbb N}}

\renewcommand{\L}{{\cal L}}
\newtheorem{thm}{Theorem}
\newtheorem{gypo}[thm]{Conjecture}

\date{}

\title{Controlling Lipschitz functions}

\author{Andrey Kupavskii\thanks{EPFL, Lausanne and MIPT, Moscow. Supported in part by the grant N 15-01-03530 of the Russian Foundation for Basic Research. E-mail: {\tt kupavskii@ya.ru}.} \and J\'anos Pach\thanks{EPFL, Lausanne and R\'enyi Institute, Budapest. Supported by Swiss National Science Foundation Grants 200020-162884 and 200021-165977. E-mail: {\tt pach@cims.nyu.edu}.} \and G\'abor Tardos\thanks{R\'enyi Institute and Central European University, Budapest. Supported by the Cryptography ``Lend\"ulet'' project of the Hungarian Academy
of Sciences and by the National Research, Development and Innovation Office, NKFIH, projects K-116769 and SNN-117879.}}
\date{}

\begin{document}
\maketitle

\noindent
\begin{abstract}
Given any positive integers $m$ and $d$, we say the a sequence of points $(x_i)_{i\in I}$ in $\R^m$ is {\em Lipschitz-$d$-controlling} if one can select suitable values $y_i\; (i\in I)$ such that for every Lipschitz function $f:\R^m\rightarrow \R^d$ there exists $i$ with $|f(x_i)-y_i|<1$.  We conjecture that for every $m\le d$, a sequence $(x_i)_{i\in I}\subset\R^m$ is $d$-controlling if and only if $$\sup_{n\in\N}\frac{|\{i\in I\, :\, |x_i|\le n\}|}{n^d}=\infty.$$ We prove that this condition is necessary and a slightly stronger one is already sufficient for the sequence to be $d$-controlling. We also prove the conjecture for $m=1$.
\end{abstract}

\section{Introduction}

The following question, in some sense dual to Tarski's famous plank problem~\cite{Ta32,McMS14,Mo32}, was raised by L\'aszl\'o Fejes T\'oth~\cite{FT74}: What is the ``sparsest'' sequence of points in the plane with the property that every straight line $\ell$ comes closer than $1$ to at least one of its points? Erd\H os and Pach~\cite{EP80} answered this question by showing that for every infinite sequence of positive numbers $(r_i)_{i\in I}$,  one can find points $p_i$ with $|p_i|=r_i$ such that every line $\ell$ passes at distance less than $1$ from $p_i$, for at least one $i\in I$, if and only if $\lim_{n\to \infty} r_n = \infty$ and $\sum_{i=1}^n\frac{1}{r_i}=\infty.$
\smallskip

Makai and Pach~\cite{MaP83} proposed a closely related, but more general question. Given a family $\cal F$ of real functions $f:\R\to\R$, we say that an infinite sequence $x_i, i\in I,$ is {\em $\cal F$-controlling} if one can choose reals $y_i, i\in I,$ such that the graph of any function $f\in{\cal F}$ ``comes close'' to at least one of the points $p_i=(x_i,y_i), i\in I,$ in the sense that
$$|f(x_i)-y_i|<1\;\;\; \mbox{\rm holds for some}\;\;\; i\in I.$$
In particular, they proved that if $\cal F$ is the family of all linear functions $f(x)=a_0+a_1x\; (a_0,a_1\in \R)$, a sequence of numbers $x_i\ge 1$ is $\cal F$-controlling if and only if $\sum_{i\in I}\frac{1}{x_i}=\infty.$ Kupavskii and Pach~\cite{KP16} managed to generalize this statement to the case where $\cal F$ consists of all polynomials $f(x)=a_0+a_1x+a_2x^2+\ldots+a_kx^k$ of degree at most $k$, for some positive $k$. In this case, the corresponding necessary and sufficient condition is $\sum_{i\in I}\frac{1}{x_i^k}=\infty.$
\smallskip

The aim of this note is to investigate the analogous problem for another interesting class of functions. Given two positive integers $m$ and $d$, let ${\cal L}(m,d)$ denote the class of {\em Lipschitz functions} from $\R^m$ to $\R^d$, that is, the class of functions for which there exists a constant $C$ such that
$$|f(x)-f(x')|\le C|x-x'|\;\;\; \mbox{\rm for all}\;\;\; x,x'\in\R^m.$$
If a function $f$ satisfies the condition above with a fixed $C>0$, then $f$ is called a {\em $C$-Lipschitz} function (or a function with {\em Lipschitz constant} $C$). Note that in this definition we can use any norm equivalent to the Euclidean norm. Throughout this note, we will work with maximum norm. For convenience, $|.|$ will stand for the {\em maximum norm}, and the word ``ball'' will refer to a ball in the maximum norm, that is, a cube.
\medskip

\noindent{\bf Definition.} {\em Given a {\em function} $f: \R^m\to \R^d$ and two points $x\in\R^m, y\in\R^d,$ we say that the pair $(x,y)$ {\em controls} $f$ if $|f(x)-y|<1$.

An infinite {\em sequence} $(x_i)_{i\in I}$ of points in $\R^m$ is said to be {\em ${\cal L}(m,d)$-controlling} or, in short, {\em $d$-controlling} if one can choose points $(y_i)_{i\in I}$ in $\R^d$ such that for every $f\in {\cal L}(m,d)$ there exists $i$ such that $(x_i,y_i)$ controls $f$.}\medskip

It follows from the definition that replacing the condition $|f(x)-y|<1$ by the inequality $|f(x)-y|<\varepsilon$ for any fixed $\varepsilon>0$, does not effect whether a sequence is $d$-controlling. To see this, it is enough to notice that $f\in{\cal L}(m,d)$ if and only if $\varepsilon f\in{\cal L}(m,d)$, and that $(x_i,y_i)$ controls $f\in{\cal L}(m,d)$ if and only if $|(\varepsilon f)(x_i)-(\varepsilon y_i)|<\varepsilon.$
\smallskip

Obviously, if a sequence is $d$-controlling, then it is also $d'$-controlling for every $1\le d'\le d$. Indeed, $\R^{d'}$ can be regarded as a subspace of $\R^d$, so every Lipschitz function from $\R^m$ to $\R^{d'}$ is a Lipschitz function from $\R^m$ to $\R^d$.
\smallskip

We solve a problem in~\cite{MaP83} by giving, for any $d$, a necessary and sufficient condition for a sequence of points in $\R$ to be $d$-controlling ($m=1$). We conjecture that this result generalizes to sequences of points in $\R^m$, for any $m\le d$, but we can prove only a slightly weaker statement.
\smallskip

The following theorem gives a necessary condition. A somewhat weaker result was
established in \cite[Theorem 3.6A]{MaP83} (it is stated in the concluding remarks of this note).

\begin{thm}\label{nec} Let $m, d$ be positive integers. If a sequence of points $(x_i)_{i\in I}$ in $\R^m$ is $d$-controlling, then we have
$$\sup_{n\in\N}\frac{|\{i\in I\, :\, |x_i|\le n\}|}{n^d}=\infty.$$
\end{thm}

Our next result shows that for $m=1$, the necessary condition in Theorem~\ref{nec} is also sufficient for a sequence of points in $\R$ to be $d$-controlling.

\begin{thm}\label{spec} Let $d$ be a positive integer. A sequence of points $(x_i)_{i\in I}$ in $\R$ is $d$-controlling if and only if
$$\sup_{n\in\N}\frac{|\{i\in I\, :\, |x_i|\le n\}|}{n^d}=\infty.$$
\end{thm}

For $m>d$, the condition in Theorems~\ref{nec} and~\ref{spec} is necessary, but {\em not sufficient} for a sequence in $\R^m$ to be $d$-controlling. To see this, observe that the sequence $(x_i)_{i\in I}$ consisting of all integer points in $\R^m$ satisfies the condition for all $d<m$. Nevertheless, this sequence is not even $1$-controlling. Indeed, for any function $h:I\to\{-1,1\}$, there exists a 2-Lipschitz function $f_h:\R^m\to\R$ for which $f(x_i)=h(i)$ for all $i\in I$. For any sequence of reals $(y_i)$, choose $\bar{h}(i)\in\{-1,1\}$ so that $|\bar{h}(i)-y_i|\ge 1$ for every $i\in I$, and notice that $f_{\bar{h}}$ is not controlled by any pair $(x_i,y_i)$.
\smallskip

However, we believe that for $m\le d$, the above condition is not only necessary but also sufficient for a sequence in $\R^m$ to be $d$-controlling.

\begin{gypo}\label{conj1} Let $m, d$ be positive integers, $m\le d$. A sequence of points $(x_i)_{i\in I}$ in $\R^m$ is $d$-controlling if and only if
$$\sup_{n\in\N}\frac{|\{i\in I\, :\, |x_i|\le n\}|}{n^d}=\infty.$$
\end{gypo}

We cannot prove this conjecture for $m>1$, but we can formulate a slightly stronger condition that is already sufficient for a sequence to be $d$-controlling, provided that $m\le d$.

\begin{thm}\label{gen} Let $m, d$ be positive integers, $m\le d$. Suppose that a sequence of points $(x_i)_{i\in I}$ in $\R^m$ satisfies the following condition for every positive $\alpha$: The set of all points $x\in\R^m$ with
$$|\{i\in I\, :\,|x_i-x|<\alpha\}|<|x|^{d-m}$$
is bounded.

Then the sequence $(x_i)_{i\in I}$ is $d$-controlling.
\end{thm}

For any $\beta>\alpha>0$, the region $\{x\in \R^m:\beta\le |x|\le 2\beta\}$ contains at least some positive constant $c=c(m)$ times $(\beta/\alpha)^m$ pairwise disjoint balls of radius $\alpha$ (that is, cubes of side length $2\alpha$, in the maximum norm). If the condition of the last theorem is satisfied, then each of these balls contains at least $\beta^{d-m}$ points $x_i$, provided that $\beta>\beta(\alpha)$ is sufficiently large. Thus, in this case,
$$|\{i\in I\, :\, |x_i|\le 2\beta\}|\ge c(\beta/\alpha)^m\beta^{d-m}=(c/\alpha^m)\beta^d.$$
Letting $\alpha\rightarrow 0,$ we obtain that the condition in
Conjecture~\ref{conj1} also holds. Roughly speaking, the condition in
Theorem~\ref{gen} is equivalent to the condition in Conjecture~\ref{conj1} for
``uniformly distributed'' sequences $x_i$, but the two conditions differ when
the density of the point sequence depends ``unevenly'' on the location. We
remark that our Theorem~\ref{nec} differs from Theorem 3.6A in \cite{MaP83} in the same sense: for ``uniform'' sequences the two
statements are equivalent, but in general they are not.
\smallskip

The exponent of $|x|$ in the right-hand side of the displayed formula in Theorem~\ref{gen} cannot be replaced by any smaller number, as follows from Theorem~\ref{nec}. Theorem~\ref{gen} disproves a conjecture from \cite{MaP83}; see the Remark at the end of Section 4.\\

It is easy to see that the sufficient condition stated in Theorem~\ref{gen} is {\em not necessary} even if $m=1$ and $d$ is arbitrary. The sequence of points consisting of $k2^{kd}$ copies of $2^k\in\R$ for every positive integer $k$, satisfies the condition of Theorem~\ref{spec} and is, therefore, $d$-controlling. On the other hand, apart from those $x\in\R$ that are closer than $\alpha$ to some power of $2$, every $x\not=0$ satisfies the inequality in Theorem~\ref{gen}. The set of these $x$ is {\em unbounded}, thus Theorem~\ref{gen} is not applicable. Since every $d$-controlling sequence of points in $\R$ can be regarded as a $d$-controlling sequence of points in $\R^m$ for any $m>1$, we obtain that the sufficient condition stated in Theorem~\ref{gen} is not necessary for a sequence to be $d$-controlling, for any values of $m$ and $d$.
\smallskip

Nevertheless, for some ``natural'' classes of sequences, the two conditions are equivalent, that is, Conjecture~\ref{conj1} holds. For instance, let $m\le d$ and $c>0$ be fixed, and consider the sequence of all points $(x_i)_{i\in I}$ in $\R^m$ whose each coordinate is the $c$-th power of some natural number. It is easy to see that this sequence satisfies both the condition in Conjecture~\ref{conj1} and the one in Theorem~\ref{gen} if $c<m/d$ and neither of them, otherwise.

\smallskip

Concerning the case $m>d$, we have a conjecture that (roughly speaking) states that a sequence in $\R^m$ is $d$-controlling if and only if there is a $d$-dimensional Lipschitz surface passing through a subset of its points that already guarantees this property. The precise statement can be formulated for every $m$ and $d$, but for $m\le d$ the conjecture is obviously true.

\begin{gypo}\label{conj2}
 Let $m, d$ be positive integers. A sequence of points $(x_i)_{i\in I}$ in $\R^m$ is $d$-controlling if and only if there exist a Lipschitz map $g:\R^d\to\R^m$ and a $d$-controlling sequence of points $(x'_i)_{i\in I'}$ in $\R^d$ with $I'\subseteq I$ such that $g(x'_i)=x_i$ for all $i\in I'$.
\end{gypo}

The ``if'' part of the conjecture is trivially true. Indeed, suppose that a sequence $(y_i)_{i\in I'}$ in $\R^d$ shows that $(x'_i)_{i\in I'}$ is $d$-controlling. Then the same sequence also shows that the sequence of points $(x_i)_{i\in I'}$ in $\R^m$ is also $d$-controlling. To see this, take any Lipschitz function $f:\R^m\to\R^d$, and observe that $f(g(x)):\R^d\to\R^d$ is also a Lipschitz function. Thus, we have $|f(x_i)-y_i|=|f(g(x'_i))-y_i|<1$ for some $i\in I'$.
\smallskip

The ``only if'' part of the conjecture evidently holds for $m\le d$. Indeed, choose $g:\R^d\to\R^m$ to be the projection to the subspace induced by the first $m$ coordinates, set $I'=I$ and $x'_i=x_i\times0^{d-m}\in\R^d$ for every $i\in I$. The important part of the conjecture is the ``only if'' direction where $m>d$.
\smallskip

The proofs of Theorems \ref{nec}, \ref{spec}, and \ref{gen}, are presented in Sections 2, 3, and 4, respectively.

\section{Proof of Theorem~\ref{nec}}

As mentioned in the introduction, a somewhat weaker statement (Theorem 3.6A) was proved in \cite{MaP83}. Here we extend the proof to the general case.

Consider a sequence $(x_i)_{i\in I}$ that violates the condition in the theorem, that is, for which
$$\sup_{n\in\N}\frac{|\{i\in I\, :\,|x_i|\le n\}|}{n^d}<\infty.$$
Given any sequence $(y_i)_{i\in I}$ of points in $\R^d$, we have to find a Lipschitz function $f\in{\cal L}(m,d)$ from $\R^m$ to $\R^d$ that is not controlled by any of the pairs $(x_i,y_i), i\in I$. We will find such a function $f$ with the property that $f(x)=g(|x|)$ for some Lipschitz function $g:\R\to\R^d$. Then it is enough to guarantee that no pair $(|x_i|,y_i)$ controls $g$. In other words, it is enough to prove the statement for $m=1$.
\smallskip

For technical reasons, we deal with the indices $i$ for which $x_i=0$, separately. Let $k$ denote the number of such indices. It follows from the assumption that $k$ is finite. Suppose without loss of generality that the index set $I$ is the set of integers larger than $-k$ and that $|x_i|$ is monotonically increasing in $i$ with $\lim_{i\rightarrow\infty}|x_i|=\infty$. Thus, we have
\[ \begin{array}{lcr}
x_i=0 & \text{if}\ i\le 0,\\
|x_i|>0 & \text{if}\ i>0.
\end{array}\]

Let  $\alpha=\sup_{i>0}\frac{i}{|x_i|^d}$ and $\beta=2\alpha^{1/d}$. For $\mu\in \R,$ denote $\lfloor \mu\rfloor$ and $\lceil \mu\rceil$ the lower and the upper integer part of $\mu$, respectively. Define a real number $\mu:=\max_{j>0}\frac{\lceil |x_j|\rceil^d}{|x_j|^d}$. It is easy to see that $\mu$ is a finite number bigger than $1$. Notice that $\alpha<\infty$ and, hence, $\beta<\infty$, because
$$\alpha=\sup_{i>0}\frac{i}{|x_i|^d}\le
\sup_{i>0}\frac{|\{j\in I\, :\,|x_j|\le |x_i|\}|}{|x_i|^d}\le \mu
\sup_{n\in\N}\frac{|\{j\in I\, :\,|x_j|\le n\}|}{n^d}<\infty.$$
\smallskip

In what follows, we define a nested sequence $\L_0\supseteq\L_1\supseteq\L_2\supseteq\ldots$ of families of $\beta$-Lipschitz functions from $\R$ to $\R^d$, we show that their intersection is nonempty, and any function $g\in\bigcap_{i\ge 0}\L_i$ meets the requirements.
\smallskip

Fix a point $y\in\R^d$ such that $|y-y_j|>1$ for every $j\le0$. Let $\L_0\subset\L(1,d)$ denote the family of all $\beta$-Lipschitz functions $g:\R\to\R^d$ with $g(0)=y$. By the choice of $y$, no function $g\in\L_0$ is controlled by any of the points $(|x_j|,y_j)$ with $j\le0$.

\begin{figure}[t!]
\begin{center}  \includegraphics[width=120mm]{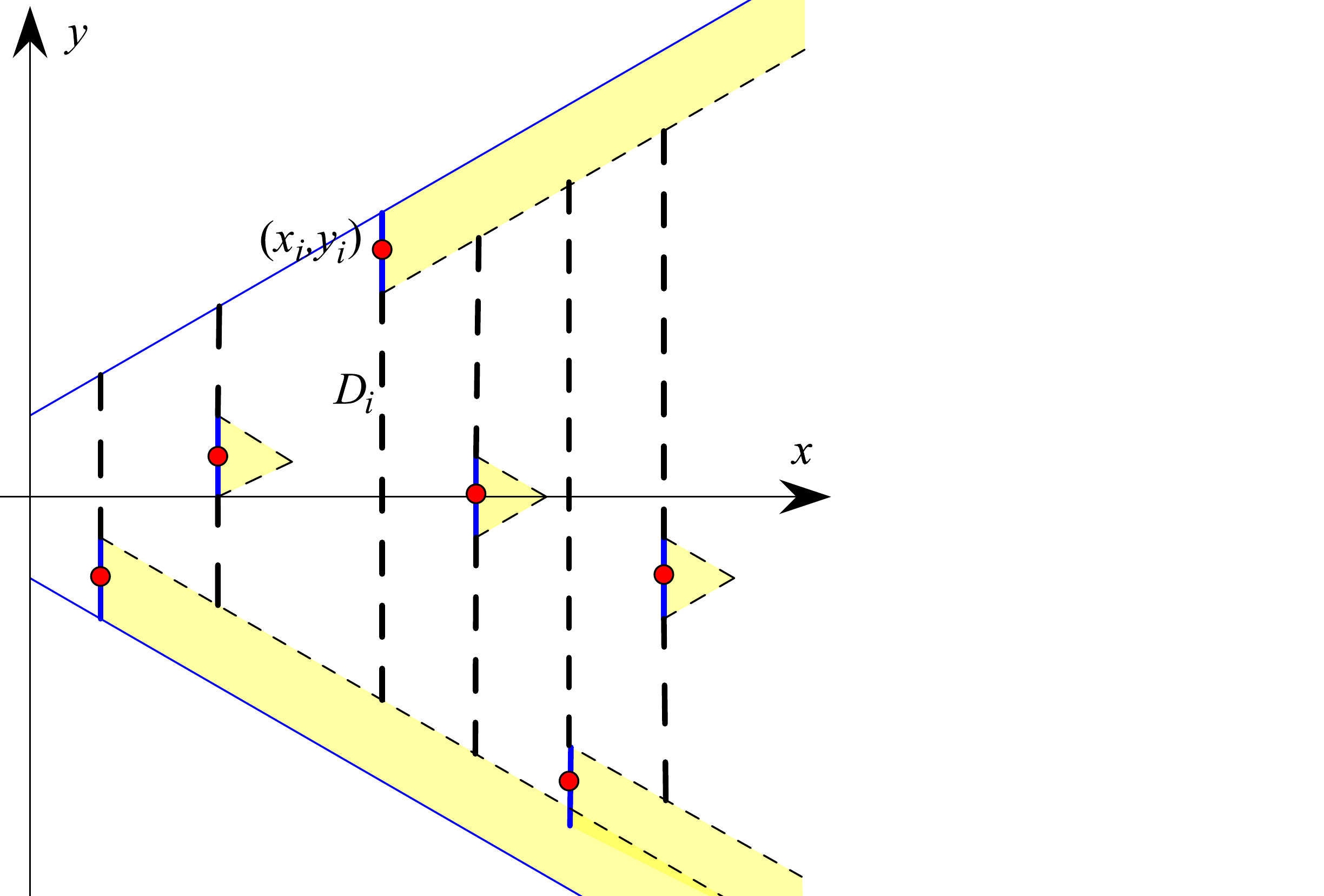}\label{pic1}  \end{center}\begin{center}  Figure 1: The case $d=1$. \end{center}
\end{figure}

For every $i>0$, let $\L_i$ be defined as the set of all functions in $\L_0$ that are not controlled by any of the pairs $(|x_j|,y_j)$ with $j\le i$, and let
$$D_i=\{g(|x_i|)\, :\, g\in\L_i\}.$$
See Fig.~1, for an illustration of the case $d=1$. (The points $(x_i,y_i)$ are marked red. If a $\beta$-Lipschitz function belongs to $\L_i$, its graph cannot intersect the yellow region incident to $(x_i,y_i)$.)

We establish a lower bound for the Lebesgue measures $\mu(D_i)$ of the sets $D_i$.

\medskip
\noindent{\bf Claim 2.1.}\; {\em For every $i\ge0$, we have\;\;
$\mu(D_i)\ge2^{d+1}\alpha|x_i|^d-2^di.$}
\medskip

\noindent{\it Proof.} By induction on $i$. For $i=0$, we have $D_0=\{y\}$, which is a nonempty set of zero measure. It follows from the definition of $\alpha$ that the bound in Claim 2.1 is strictly positive for every $i>0$. Assume that we have already verified the Claim for some $i\ge0$, and we want to prove it for $i+1$.

Let $D'=\{g(|x_{i+1}|)\, :\, g\in\L_i\}$. Clearly, $D'$ can be obtained as the Minkowski sum of $D_i$ and the ball $B_r=B_r(0)$ of radius $r=\beta(|x_{i+1}|-|x_i|)$ around the origin. On the other hand, we have $D_{i+1}=D'\setminus B_1(y_{i+1})$, where $B_1(y_{i+1})$ denotes the ball of radius $1$ around $y_{i+1}$. Therefore,
$$\mu(D_{i+1})\ge\mu(D')-\mu(B_1(y_{i+1}))=\mu(D'+B_r)-\mu(B_1(y_{i+1})).$$
By the Brunn-Minkowski inequality, we have
$$\mu(D'+B_r)\ge(\mu^{1/d}(D_i)+\mu^{1/d}(B_r))^d.$$
Combining the last two inequalities,
$$\mu(D_{i+1})\ge(\mu^{1/d}(D_i)+\mu^{1/d}(B_r))^d-\mu(B_1(y_{i+1})).$$
As we use the maximum norm, we have $\mu(B_1(y_{i+1}))=2^d$ and $\mu(B_r)=2^dr^d$. Using the inductive hypothesis, we get the following chain of implications.
\begin{alignat*}{2}
  \phantom{\geq}  \mu(D_{i+1}) & \ge 2^{d+1}\alpha|x_{i+1}|^d-2^d(i+1) &&\Leftarrow \\
    \big((2^{d+1}\alpha|x_i|^d-2^di)^{\frac 1d}+2r\big)^d-2^d &\ge 2^{d+1}\alpha|x_{i+1}|^d-2^d(i+1)\ \ \ &&\Leftrightarrow \\
        \big((2\alpha|x_i|^d-i)^{\frac 1d}+r\big)^d &\ge 2\alpha|x_{i+1}|^d-i\ \ \ &&\Leftrightarrow \\
   \beta(|x_{i+1}|-|x_i|) &\ge (2\alpha|x_{i+1}|^d-i)^{\frac 1d}-(2\alpha|x_i|^d-i)^{\frac 1d}, &&
\end{alignat*}
where $\beta=2\alpha^{1/d},$ as before. By the definition of $\alpha$, we have $2\alpha x^d-i\ge \alpha x^d$ for every $x\ge |x_i|$. Consider the function $f(x):=(2\alpha x^d-i)^{1/d}$. Then
$$f'(x)= \frac 1d 2\alpha d\frac{x^{d-1}}{(2\alpha x^d-i)^{\frac{d-1}d}}\le 2\alpha \frac{x^{d-1}}{(\alpha x^d)^{\frac{d-1}d}}=2\alpha^{1/d},$$
for every $x\ge |x_i|$.
Therefore, the last inequality of the chain holds, and so does the first one, as claimed. \;Q.E.D.
\medskip

In particular, it follows from Claim 2.1 that $D_i\not=\emptyset$ and, hence, $\L_i$ is not empty for every $i\ge0$. To complete the proof of Theorem 1, it is enough to note that the set $\L_0$ is {\em compact} in the pointwise topology. Therefore, $\bigcap_{i\ge 0}\L_i\neq\emptyset$. By definition, no function $g\in \bigcap_{i\ge 0}\L_i$ is controlled by any pair $(|x_i|,y_i)$, as required.

\section{Proof of Theorem~\ref{spec}}

The ``only if'' part of the theorem is a special case of Theorem~\ref{nec}. Thus, we have to prove only the ``if'' part.
\smallskip

Let $(x_i)_{i\in\N}$ be a sequence of real numbers satisfying the ``density condition''
$$\sup_{n\in\N}\frac{|\{i\in I\, :\,|x_i|\le n\}|}{n^d}=\infty.$$
Split this sequence into two sequences, one consisting of the nonnegative numbers and the other consisting of the negative ones. At least one of these two sequences must satisfy the above density condition, so we can assume without loss of generality that, say, $x_i\ge0$ for all $i$.
If $(x_i)_{i\in\N}$ has a convergent subsequence $(x_{i_j})_{j\in \N}\rightarrow x$, as $j\rightarrow\infty$, then choose any sequence of points $(y_j)_{j\in\N}$, everywhere dense in $\R^d$. Obviously, every Lipschitz function $f:\R\rightarrow\R^d$ is controlled by infinitely many pairs $(x_{i_j},y_j)$. Therefore, we can assume without loss of generality that $(x_i)_{i\in\N}$ is an increasing sequence of nonnegative numbers, tending to infinity.
\smallskip

We need a simple statement about a finite portion of the sequence $(x_i)$.
\smallskip

Fix a positive integer $j$. Let $\L_j$ denote the family of $j$-Lipschitz functions $f:\R\to\R^d$ with $|f(0)|\le j$. (Note that this deviates from the definition of $\L_j$ used in the previous section.) We also fix $n\in\N$. Let $k=k(j,n)=(j(n+1)+1)^d$, and assume that $x_i\le n$ for $1\le i\le k$.
\smallskip

Since we use the maximum norm in $\R^d$, the ball (cube) $B_r$ of radius $r=j(n+1)$ around the origin can be uniquely partitioned into $k$ balls of radius $r'=\frac{r}{r+1}<1$. Let the centers of these balls be denoted by $z_i$ and the balls themselves by $B_{r'}(z_i),\;1\le i\le k.$ Index the centers $z_i$ decreasingly with respect to the lexicographic order. For every $i, 1\le i\le k$, set $$y_i=z_i-j(n-x_i)v,$$
where $v$ is the all-$1$ vector in $\R^d$. See Fig. 2, which depicts the case $d=1$.

\begin{figure}[t!]
\begin{center}  \includegraphics[width=120mm]{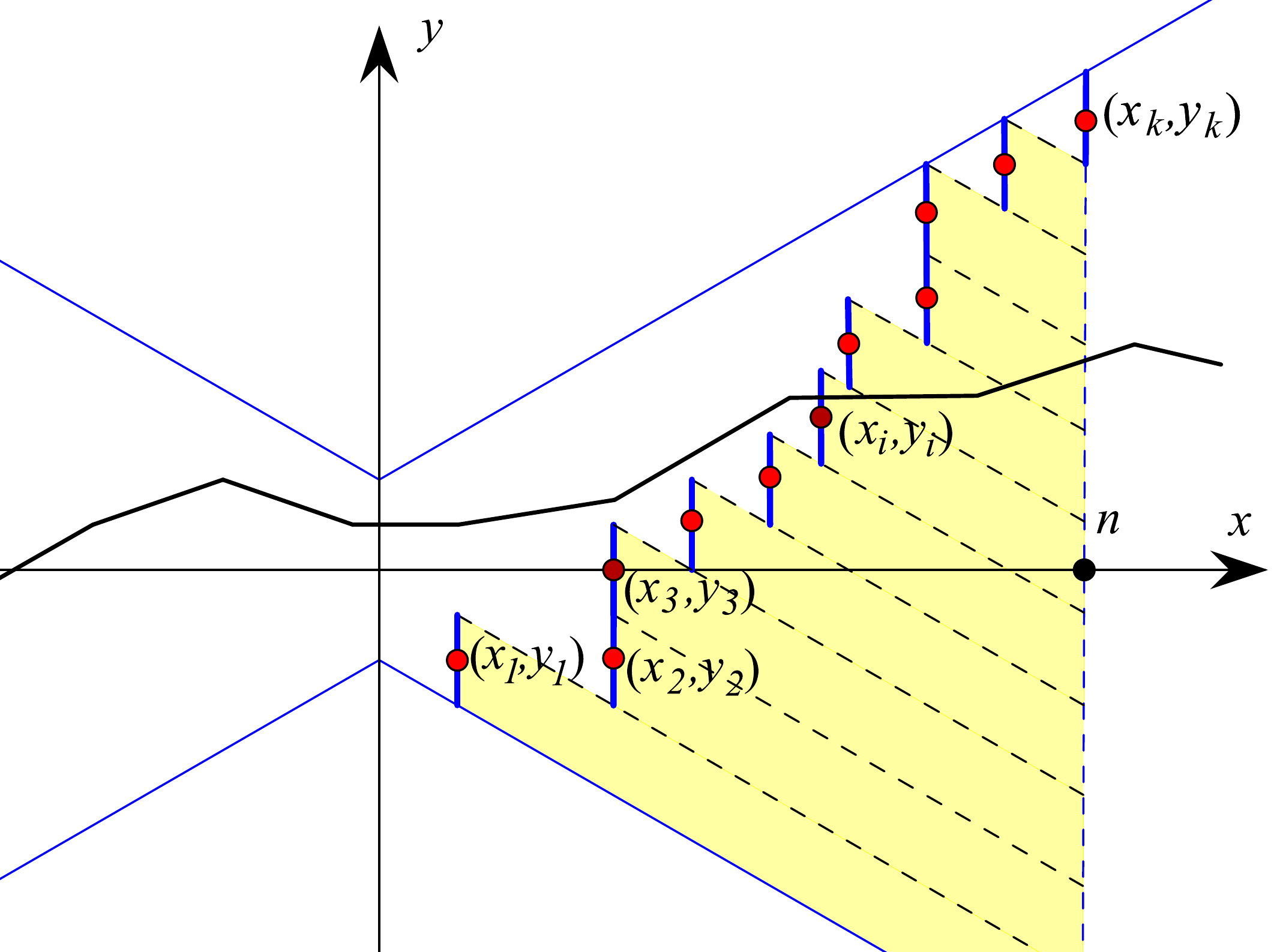}\label{pic2}  \end{center}\begin{center}  Figure 2. \end{center}
\end{figure}
\medskip

\noindent{\bf Claim 3.1.} {\em Any function $f\in\L_j$ is controlled by one of the pairs $(x_i,y_i), \; 1\le i\le k.$}
\medskip

\noindent{\em Proof.} Let $f$ be an arbitrary element of $\L_j$. Notice that the function $g(x)=f(x)+j(n-x)v$ is monotonically decreasing in all of its coordinates, and that $g(x)\in B_r$ for every $x\in [0,n].$

Consider the set $S$ of all indices $1\le i\le k$ such that $g(x_i)$ is contained in a ball $B_{r'}(z_{i'})$ for some $1\le i'\le i$. As $g(x_k)$ is in $B_r$, it belongs to a ball $B_{r'}(z_{i'})$ for some $1\le i'\le k$. Therefore, $k\in S$, so that the set $S$ is not empty. Let $i_0$ denote the smallest element of $S$.

Then we have $g(x_{i_0})\in B_{r'}(z_{i_0})$. Indeed, otherwise $g(x_{i_0})\in B_{r'}(z_{i'})$ for some index $i'<i_0$. Thus, $i_0>1$. Using the monotonicity of $g$ and the monotonicity of the sequences $(x_i)_{1\le i\le k}$ and $(z_i)_{1\le i\le k}$, we obtain that $g(x_{i_0-1})\in B_{r'}(z_{i''})$ for some $i''\le i'\le i_0-1$, contradicting the minimality of $i_0$. Hence, $$1>r'\ge|g(x_{i_0})-z_{i_0}|=|f(x_{i_0})-y_{i_0}|.$$
This means that $(x_{i_0},y_{i_0})$ controls $f$, as claimed. \; Q.E.D.

\medskip

Now we can easily finish the proof of Theorem~\ref{spec}. We need to show that the sequence $(x_i)_{i\in\N}$ is $d$-controlling. To control all functions in $\L_j$ for a fixed $j$, pick an $n=n(j)$ such that for at least $k=(j(n+1)+1)^d$ distinct indices $i$ we have $x_i\le n$. It follows from the density condition that such an $n$ exists.
\smallskip

By Claim 3.1, we can choose $y_i$ for $k$ distinct indices $i$ such that every function in $\L_j$ is controlled by one of the $k$ pairs $(x_i,y_i)$. Repeat this step this successively for $j=1,2,\ldots$, making sure that we always use pairwise disjoint sets of indices. This is possible, because removing any finite number of elements from $(x_i)_{i\in \N}$, the remaining sequence still satisfies the density condition. Since every Lipschitz function $\R\to\R^d$ belongs to one of the classes $\L_j$, after completing the above process for all $j\in\N$, all Lipschitz functions $\R\to\R^d$ will be controlled by one of the pairs $(x_i,y_i)$. This proves Theorem~\ref{spec}.

\section{Proof of Theorem~\ref{gen}}

As in the proof of Theorem \ref{spec}, for every positive integer $j$, $\L_j$
denotes the family of $j$-Lipschitz functions $f:\R^m\to\R^d$ with $|f(0)|\le
j$. As we did in that proof, we fix $j$ and we show that one can control $\L_j$ using only finitely many points $x_i$. To complete the proof of Theorem~\ref{gen}, we perform this step for $j=1,2,\ldots$, sequentially, observing that the density condition in the theorem continues to hold even if we delete any finite number of points $x_i$ from our sequence.

%
%
%
%

\medskip
The proof of Theorem~\ref{gen} is based on a topological lemma. We consider a continuously moving set $D$ that leaves a ball $B\subset R^d$. By continuity, each point of $D$ must cross the boundary of the ball. Using Brouwer's fixed point theorem, we find a point $z\in D$ that crosses the boundary at a point with a special property. See Figure 3, for an illustration. The color gradation distinguishes different points of $D$, that is, points of the same color indicate the trajectory of a point, as it progresses in time $t$.
\medskip

\begin{figure}[t!]
\begin{center}  \includegraphics[width=160mm]{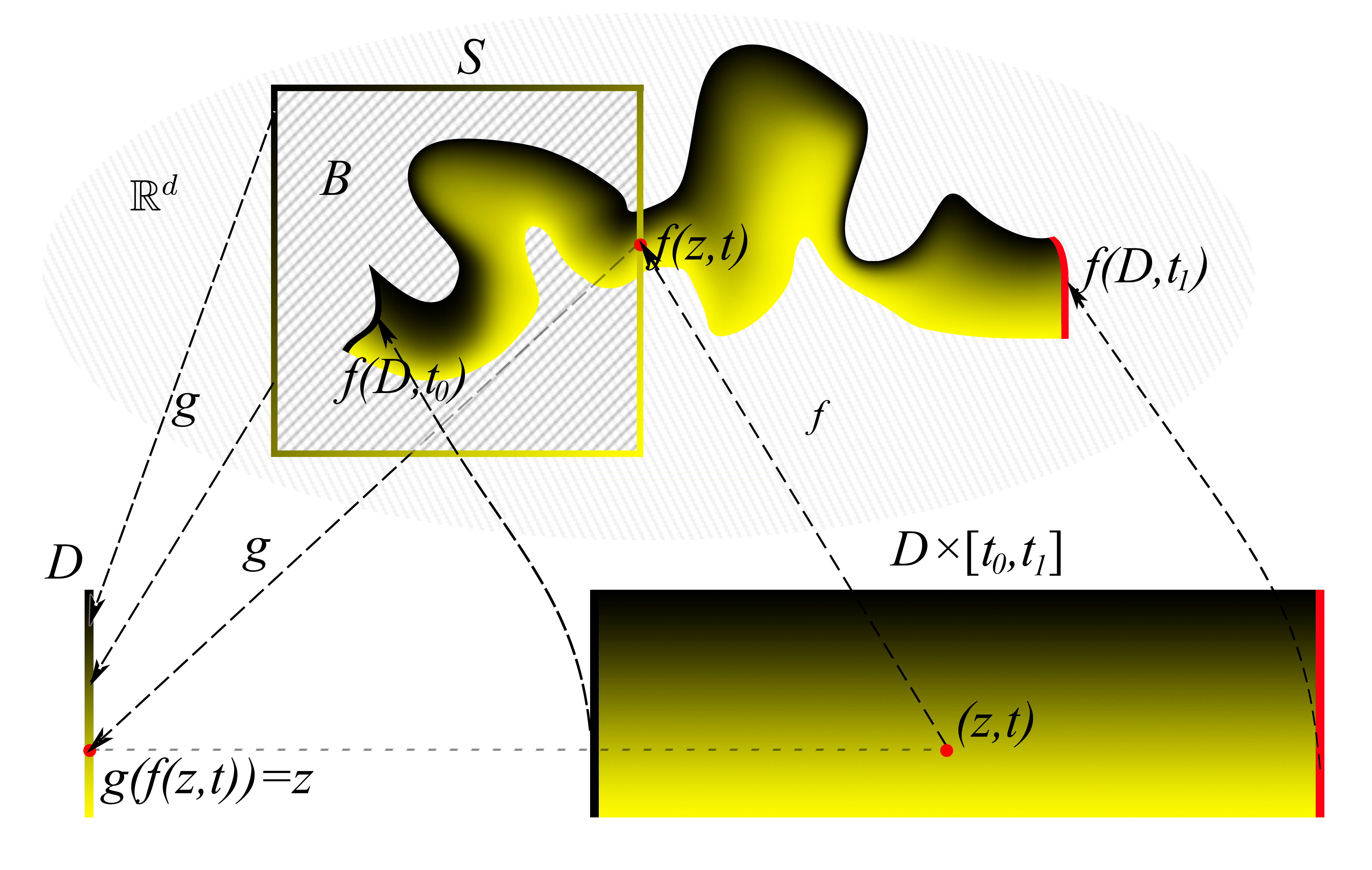}\label{pic3}  \end{center}\begin{center}  Figure 3. \end{center}
\end{figure}
\medskip

\noindent{\bf Lemma 4.1.} {\em  Let $d$ be a positive integer. Let $B$ denote a closed ball of positive radius around the origin in $\R^d$, and let $S$ stand for the boundary of $B$. Let $J=[t_0,t_1]$ be a closed interval on the real line, let $D$ be an arbitrary topological space, and let $f:D\times J\to\R^d$ and $g:B\to D$ be continuous functions.

If $f(z,t_0)\in B\setminus S$ and $f(z,t_1)\notin B$ for all $z\in D$, then there exist $z\in D$ and $t\in (t_0,t_1)$ such that $f(z,t)\in S$ and $g(f(z,t))=z$.}
\medskip

\noindent{\em Proof.} Denote $l$  the radius of $B$. Let $B'=B\times J\subset\R^{d+1}$, and let $h:B'\to\R^{d+1}$ be defined as $h(y,t)=(y',t')$, where $y'=f(g(y),t)$ and $t'=t-|y'|+l$.  Let $c:\R^{d+1}\to B'$ be a coordinate-wise retraction; to be specific, let $c(y,t)=(\min(1,l/|y|)\cdot y,\min(t_1,\max(t_0,t)))$. Finally, let $\bar h:B'\to B'$ be the composition of these functions: $\bar h(y,t)=c(h(y,t))$.

Clearly, $\bar h$ is continuous and $B'$ is homeomorphic to the $(d+1)$-dimensional ball. Thus, we can apply Brouwer's fixed point theorem to conclude that there exists $(y,t)\in B'$ with $\bar h(y,t)=(y,t)$. Let $z=g(y)$ and $(y',t')=h(y,t)$.

If $t=t_0$, then the second coordinate of $c(y',t')$ equals $t_0$. Hence, either we have $t'=t_0$ or $t'$ was retracted to $t_0$ from the left. Since $t'\le t$, we have $|y'|\ge l$, and thus $f(z,t_0)\notin B\setminus S$, contradicting our assumption.

Analogously, if $t=t_1$, then $t'\ge t$, so $|y'|\le l$, implying that $f(z,t_1)\in B$, which is again  a contradiction.

Consequently, we must have $t_0<t<t_1$. Using the fact that $c$ is a retraction, we obtain that $t'=t$, so $|y'|=l$ and $y'=f(z,t)\in S$. We must also have $y'=y$, which implies that $g(f(z,t))=g(y)=z$. \; Q.E.D.
\medskip

To apply the lemma, we think of $\R^m$ as a product space $\R^{m-1}\times \R$,
with the last coordinate considered as time.
Let $D\subset \R^{m-1}$ and  $B\subset \R^d$ be balls, let $J=[t_0,t_1]$ be an interval, and let $g:B\to D$ a linear map (see Fig. 3). Consider any $j$-Lipschitz function $f:\R^m\to\R^d\; (m\le d),$ and focus our attention on the restriction of $f$ to $D\times J$. In order to apply Lemma~4.1, we choose $B$ large enough to make sure that $f(z,t_0)$ lies in the interior of $B$ for all $z\in D$. By the lemma, we can either find $z\in D$ such that $x=(z,t_1)$ satisfies $y=f(x)\in B$, or there exists $x=(z,t)\in D\times J$ such that $y=f(x)$ belongs to the boundary of $B$ and $g(y)=z$. Our goal is to find sufficiently many indices $i\in I$ with $x_i\in D\times J$, and to assign appropriate values $y_i$ to them, so that for every conceivable pair $(x,y)$ provided by the lemma we can find a pair $(x_i,y_i)$ that is close to it. Specifically, if we have $|x-x_i|<\frac 1{2j}$ and $|y-y_i|<\frac 12$ for some $i$, then $f(x)=y$ implies that the pair $(x_i,y_i)$ controls $f\in\L_j$. Next we spell out the details of proof.
\medskip

\noindent
\textit{Proof of Theorem~\ref{gen}. \ } Let $m\le d$ and let $(x_i)_{i\in I}$
be a sequence of points in $\R^m$ satisfying the density condition in the
theorem. Let us fix $j\in \N$. As we have pointed out earlier, it is
sufficient to show that we can select finitely many indices $i\in I$ and
assign to them suitable points $y_i\in\R^d$ such that every function in $\L_j$
is controlled by at least one of the pairs $(x_i,y_i)$.
\smallskip

Set $\epsilon=1/(8j+8)$ and choose a positive integer $c$ with
$c^m>4d/\epsilon^{d-m}$. Using the density condition in the theorem with $\alpha=\epsilon/c$,
we obtain that there exists $t_0>j+1$ such that $$|\{i\in
I\,:\,|x_i-x|<\alpha\}\ge|x|^{d-m}$$
holds for every $x\in \R^m$ with $|x|\ge t_0-2\epsilon$.  Set $l=\lfloor jt_0+j\rfloor+1<(j+1)t_0$. The density
condition we really need for our argument is
$$|\{i\in I\,:\,|x_i-x|<\epsilon\}\ge4d(2l)^{d-m},$$
which holds for every $|x|\ge t_0-\epsilon$, since the ball of radius $\epsilon$
around $x$ can be split into $c^m$ internally disjoint balls of radius
$\alpha$, each containing at least $(|x|-\epsilon)^{d-m}$ points $x_i$ in their
interior. This adds up to total of $c^m(t_0-2\epsilon)^{d-m}>4d(2l)^{d-m}$, as
required. Finally, set $t_1>t_0$ such that
$$|\{i\in I\,:\,|x_i-x|<\epsilon\}|\ge4d(2l)^{d-m}+(2l)^d$$
holds for all $|x|\ge t_1-\epsilon$. The existence of such a value $t_1$ follows easily from the density condition on the sequence $(x_i)$, because for every sufficiently large $|x|$, the
left-hand side of the above inequality is at least $|x|^{d-m}$, while its right-hand side is a
constant.
\smallskip

Let $D=\{z\in\R^{m-1}\, :\,|z|\le t_0\}$ be the ball of radius $t_0$ around the
origin in $\R^{m-1}$, let $J=[t_0,t_1]$, let $B=\{y\in\R^d\, :\,|y|\le l\}$ be
the ball of radius $l$ around the origin in $\R^d$, and let $S=\{y\in\R^d\,
:\,|y|=l\}$ denote the sphere bounding $B$. Define a linear map $g:B\to D$ by
setting
$$g(y_1,\ldots,y_d)=\frac {t_0}{2l}(y_1-y_m,y_2-y_m,\ldots,y_{m-1}-y_m).$$
We identify $\R^m$ with $\R^{m-1}\times\R$ and will use the notation $(z,t)\in\R^m$ for $z\in\R^{m-1}$ and $t\in\R$.

Cover $D\times J$ with internally disjoint balls (cubes, in the $l_{\infty}$-norm) of radius $\epsilon$. These balls will be referred to as the
{\em $\epsilon$-balls}. Let $Z=Z_0\times Z_1$ be a fixed $\epsilon$-ball, where
$Z_0\subset\R^{m-1}$ is a ball of radius $\epsilon$ and $Z_1$ is an interval
of length $2\epsilon$.
\smallskip

The sphere $S$ consists of $2d$ facets ($d-1$-dimensional cubes). A facet
is obtained by fixing one of the $d$ coordinates to $l$ or $-l$, and letting
the other coordinates take arbitrary values in the interval $[-l,l]$. Consider all
points $y$ on a facet such that $g(y)\in Z_0$. If the fixed
coordinate of the facet is one of the first $m$ coordinates, and such points $y$
exist at all, then the first $m$ of their coordinates are determined within an
interval of $8l\epsilon/t_0\le1$, while the remaining coordinates can take
arbitrary values in $[-l,l]$. This set can be covered by at most $(2l)^{d-m}$
balls of radius $1/2$. We refer to these balls as the {\em $1/2$-balls for $Z$}.

Next, consider
all points $y$ on a facet of $S$ such that $g(y)\in Z_0$, but assume that the fixed
coordinate of this facet is one of the last $d-m$ coordinates. Cover
this set with balls of radius $1/2$, as follows. Partition the possible
values of the $m$'th coordinate into $4l$ intervals, each of length $1/2$. These
intervals determine each of the first $m-1$ coordinates of $y$ within an
interval of length $1$, and there are $d-m-1$ further coordinates that can take any value in
$[-l,l]$. We have $2(2l)^{d-m}$ balls of radius $1/2$ that cover all points
$y$ on this facet with $g(y)\in Z_0$. Summing up over all facets of $S$, we have at most $4d(2l)^{d-m}$ $1/2$-balls for $Z$. For each of these
$1/2$-balls $W$ for $Z$, select a separate index $i\in I$ such that $x_i$ lies in the
interior of $Z$, and set $y_i$ to be the center of the ball $W$. Note that the center $x$ of $Z$ satisfies $|x|\ge t_0-\epsilon$ (otherwise, $Z$ would be disjoint from $D\times J$). Thus, by our choice of $t_0$, we have enough indices to choose from. We repeat the same procedure for every the $\epsilon$-ball $Z$.

\smallskip

\noindent{\em Case 1:} Consider now any $f\in{\cal L}_j$ for which there exists $x=(z,t)\in D\times J$ such that $y=f(x)\in S$ and $g(y)=z$.

Clearly, $x\in Z$ for some $\epsilon$-ball $Z$, and $y\in W$ for some $1/2$-ball $W$ for $Z$. Thus, there exists $i\in I$ such that $x_i$ lies in the interior of $Z$, and $y_i$ is the center of $W$. This implies that $|x_i-x|<2\epsilon$ and $|y_i-y|\le\frac 12$. Using the Lipschitz property, we obtain
$$|f(x_i)-y|=|f(x_i)-f(x)|\le j|x_i-x|<2j\epsilon<\frac12.$$
Hence, $(x_i,y_i)$ controls $f$, as $|f(x_i)-y_i|\le|f(x_i)-y|+|y_i-y|<1$.

\smallskip

\noindent{\em Case 2:} It remains to deal with the case where for some $f\in{\cal L}_j$ we cannot find $x=(z,t)\in D\times J$ such that $y=f(x)\in S$ and $g(y)=z$.

Let $f$ be such a function. Notice that $f(z,t_0)\in B\setminus S$ for every $z\in D$. Indeed, we have $|(z,t_0)|=t_0$ and, hence, $|f(z,t_0)|\le jt_0+|f(0)|\le jt_0+j<l$, as required.
Then, according to Lemma~4.1, if we cannot find $x=(z,t)\in D\times J$ such that $y=f(x)\in S$ and $g(y)=z$, then $f(z,t_1)\in B$ must hold for some $z\in D$. We show that in this case one can select a few more indices $i\in I$ and set the corresponding values $y_i$ so that for some of the newly selected indices $i$, the pairs $(x_i,y_i)$ control $f$.

To achieve this, cover the entire ball $B$ with $(2l)^d$ balls of radius $1/2$, and refer to them as {\em new balls}. For any $\epsilon$-ball $Z$ that contains a point $(z,t_1)$ and for any new ball $W$, choose a separate (yet unselected) index $i\in I$ such that $x_i$ lies in the interior of $Z$, and set $y_i$ to be the center of $W$. Note that the center $x$ of $Z$ satisfies the inequality $|x|\ge t_1-\epsilon$. Thus, by our choice of $t_1$, we have enough indices to choose from. It can be shown by a simple computation similar to the above one that if for some $z\in D$ we have $y=f(z,t_1)\in B$, then for the indices $i\in I$ selected for the $\epsilon$-ball containing $(z,t_1)$ and the new ball containing $y$ the pair $(x_i,y_i)$ controls $f$.

\smallskip

This completes the proof of the fact that every $f\in{\cal L}_j$ is controlled by one of the
pairs $(x_i,y_i)$ and, hence, the proof of Theorem~\ref{gen}.  \; Q.E.D.

\bigskip
\noindent{\bf Remark.}  Makai and Pach \cite{MaP83} proved the following result (that also follows from our Theorem~\ref{nec}): Let $m\le d$ and let $A$ be a set of points in $\R^m$ satisfying the condition that for any $x\in \R^m$, the number of points in the unit ball around $x$ is at most $K(|x|^{d-m}+1)$, where $K$ is a suitable constant. Then $A$ is {\em not} $d$-controlling. Makai and Pach made the conjecture that the same statement remains valid if for any $x\in \R^m$, the unit ball around $x$ contains at most $K(|x|^{d-1}+1)$ points of $A$. This would be a significant improvement for $m>1$. However, our Theorem~\ref{gen} shows that no such improvement is possible. Indeed, for any function $f:\R^+\to\R^+$ tending to infinity, one can construct a set of point in $\R^m$ with at most $f(|x|)|x|^{d-m}$ points in the unit ball around any point $x$, but still satisfying the condition of Theorem~\ref{gen}. By the theorem, such a set is $d$-controlling.

We close this paper by constructing an explicit set $A$ with the properties mentioned above. We choose an increasing sequence of reals $c_i>4$ such that $f(x)>2^{m(i+2)+d}$ whenever $x\ge c_i-2$. Consider the set $S_{i}:=\{x\in2^{-i} \mathbb Z^m\mid c_i\le|x|<c_{i+1}\}$. Form a set $A_i\subset \R^m$ by collecting $\left\lceil|x|^{d-m}\right\rceil$ points from the ball of radius $2^{-i}$ around every point $x\in S_i$. Consider the set $A:=\cup_{i=1}^{\infty}A_i$. For $|x|>c_{i}$, we have at least $|x|^{d-m}$ points of $A$ in the $2^{1-i}$-ball around $x$. This shows that $A$ satisfies the conditions of Theorem~\ref{gen} and is, therefore, $d$-controlling. On the other hand, if $i$ is the highest index such that the unit ball around $x$ contains a point in $A_i$, then $|x|>c_i-2>2$, and the unit ball around $x$ contains at most $\left\lceil(|x|+2)^{d-m}\right\rceil$ points of $A$ around each of the at most $2^{m(i+2)}$ points of $\cup_{i=1}^\infty S_i$ in the ball of radius $2$ about $x$. By our choice of $c_i$, this shows that the unit ball around $x$ contains at most $f(|x|)|x|^{d-m}$ points of $A$, as claimed.\\

\textsc{Acknowledgements. } We thank the referee for carefully reading the text and pointing out an inaccuracy in the proof of Claim~2.1.

\end{document}